\documentclass[12pt]{article}

\def \figpath{figures}
\def \kws{Adaptive control, Robust adaptive control, Robust control, H-infinity control, Game theory}
\newcommand{\forArxiv}[1]{#1}  
\newcommand{\forJournal}[1]{}      

\forJournal{
\RequirePackage{xcolor}
	\documentclass[journal,twoside,web, final]{ieeecolor}
	\makeatletter
	\let\NAT@parse\undefined
	\makeatother
	\usepackage{etoolbox}
	\@ifundefined{color@begingroup}%
	  {\let\color@begingroup\relax
	   \let\color@endgroup\relax}{}%
	\def\fix@ieeecolor@hbox#1{%
	  \hbox{\color@begingroup#1\color@endgroup}}
	\patchcmd\@makecaption{\hbox}{\fix@ieeecolor@hbox}{}{\FAILED}
	\patchcmd\@makecaption{\hbox}{\fix@ieeecolor@hbox}{}{\FAILED}

	\markboth{\hskip25pc IEEE TRANSACTIONS AND JOURNALS TEMPLATE}
	{Kjellqvist \MakeLowercase{\textit{et al.}}: Output Feedback Minimax Adaptive Control}
	}

\forJournal{

}
\forJournal{
\usepackage{generic}
}
\usepackage{cite}
\usepackage{lipsum}
\usepackage{amsmath,amssymb,amsfonts, amsthm}
\usepackage{textcomp}
\usepackage{verbatim} 
\usepackage{tikz}
\usepackage{hyperref}
\hypersetup{
  colorlinks=true,
  linkcolor=blue,
  citecolor=blue,
  filecolor=blue,
  urlcolor=blue,
  pdftitle={Minimax Adaptive Control and Estimation},
  pdfauthor={Olle Kjellqvist},
  pdfsubject={PhD Thesis},
}
\def\BibTeX{{\rm B\kern-.05em{\sc i\kern-.025em b}\kern-.08em
    T\kern-.1667em\lower.7ex\hbox{E}\kern-.125emX}}

\usetikzlibrary{arrows,shapes,positioning,calc}

 \newcommand*{\defeq}{\stackrel{\mbox{\small{\ensuremath{\mathsf{def}}}}}{=}} 
\newcommand{\tran}{{\mkern-1.5mu\mathsf{T}}}

\newcommand{\modelset}{\mathcal{M}}

\newcommand{\Hinf}{\mathcal{H}_\infty}

\renewcommand{\S}{\mathbb{S}}

\newcommand{\R}{\mathbb{R}}

\newcommand{\seq}[3]{\mathbf{#1}_{#2}^{#3}}

\newcommand{\bmat}[1]{\begin{bmatrix} #1 \end{bmatrix}}

\DeclareMathOperator{\bdiag}{BlockDiag}

\DeclareMathOperator*{\argmax}{arg\,max}
\DeclareMathOperator{\quadform}{\sigma}

\newcommand{\B}{\mathcal{B}}

\theoremstyle{plain}

\newtheorem{theorem}{Theorem}
\newtheorem{prop}{Proposition}
\newtheorem{lemma}{Lemma}
\newtheorem{corollary}{Corollary}[theorem]
\newtheorem{remark}{Remark}

\newtheorem{assumption}{Assumption}
\newtheorem{problem}{Problem}

\begin{document}
\title{Output Feedback Minimax Adaptive Control}
\author{Olle Kjellqvist and Anders Rantzer%
\thanks{This project has received funding from the European Research Council (ERC) under the European Union's Horizon 2020 research and innovation programme under grant agreement No 834142 (ScalableControl).
Both authors are  with the Department of Automatic Control, 
Lund University, Lund, Sweden %
(e-mail: [olle.kjellqvist, anders.rantzer]@control.lth.se).
}
}
\maketitle

\begin{abstract}
	This paper formulates adaptive controller design as a minimax dual control problem. 
	The objective is to design a controller that minimizes the worst-case performance over a set of uncertain systems. 
	The uncertainty is described by a set of linear time-invariant systems with unknown parameters. 
	The main contribution is a common framework for both state feedback and output feedback control. 
	We show that for finite uncertainty sets, the minimax dual control problem admits a finite-dimensional information state. 
	This information state can be used to design adaptive controllers that ensure that the closed-loop has finite gain.
	The controllers are derived from a set of Bellman inequalities that are amenable to numerical solutions. 
	The proposed framework is illustrated on a challenging numerical example.
\end{abstract}
\forJournal{
}
\forArxiv{
}
\section{Introduction}
This paper addresses the design of adaptive controllers with guaranteed performance for linear time-invariant systems with uncertain parameters.
The performance index is quadratic with a soft constraint on the size of the disturbance.
The performance measure quantifies transient behavior and, if finite, guarantees a bounded $\ell_2$-gain.
Via the small-gain theorem, we guarantee stability in the presence of unmodeled dynamics.

\begin{figure}
	\centering
	\includegraphics{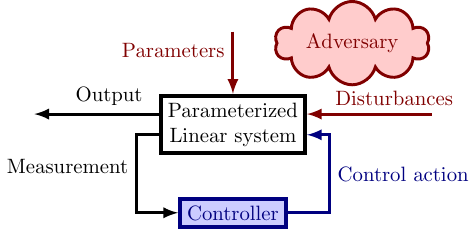}
	\caption{The minimax control problem. The controller minimizes a performance index by selecting inputs, while the adversary selects the parameter realization and disturbances to maximize it.}
	\label{fig:intro:problem1}
\end{figure}

This property implies that the closed-loop system behaves well even if the assumptions on the model class are violated, as long as the violation is minor.
Toward this end, we do not make any assumptions on the statistical properties of the parameters or the exogenous signals.
Instead, the underlying models are deterministic, and the uncertain parameters and signals are chosen by an adversary that seeks to maximize the performance index.
See Fig.~\ref{fig:intro:problem1} for an illustration of the problem.

\subsection{Contributions}
To address the complexities and challenges outlined above, this paper makes several key contributions to the field of adaptive control.
\subsubsection{Unifying State-Feedback and Output-Feedback}
\begin{figure}
	\centering
	\includegraphics{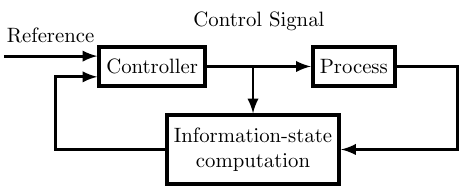}
	\caption{Closed-loop with information state-feedback control.}
	\label{fig:intro:information_state_feedback}
\end{figure}

We show that the state-feedback and output-feedback minimax dual control problems can be reduced to a minimax control problem with linear (known) dynamics and uncertain objective functions.
This reduction is based on the concept of information-state feedback\cite{Bertsekas1973Sufficient, James95Robust}, and is illustrated in Fig.~\ref{fig:intro:information_state_feedback}.
The problem with uncertain objective functions is introduced in Section~\ref{sec:principal}.
State-feedback and output-feedback minimax dual control problems and their reductions are presented in Sections~\ref{sec:state_feedback} and~\ref{sec:output_feedback}, respectively.

\subsubsection{Finite-Dimensional Information State}
We show that if the uncertain parameters belong to a finite set, the optimal output-feedback controller is observer-based and can be computed by dynamic programming.
This is a specialization of the result in~\cite{James95Robust} for the nonlinear $\Hinf$-control problem.
However, in contrast to the general result~\cite{James95Robust}, we show in Section~\ref{sec:output_feedback} that the observer state is finite-dimensional, and we provide a constructive method to compute the observer state.

\subsubsection{Heuristics for Suboptimal Controller Synthesis}
We provide a heuristic method to synthesize suboptimal controllers.
The method is based on approximating the value function as a piecewise-quadratic function and the controllers as certainty-equivalence controllers.
We introduce periodic Bellman inequalities to deal with delays and other nonminimum-phase behavior.
The method generalizes \cite[Theorem 3]{Rantzer2021L4DC} to minimax control of linear time-invariant systems with unknown objective functions.
Hence, the method is applicable to both state-feedback and output-feedback minimax dual control problems.
The resulting controllers have guaranteed worst-case performance in the sense of a bounded $\ell_2$-gain.
The details of the method are presented in Section~\ref{sec:explicit}, and the use of periodicity to deal with nonminimum-phase behavior is examplified in Section~\ref{sec:examples:statefeedback}.

\subsubsection{Numerical Examples}
We provide a \texttt{Julia}~\cite{Bezanson2017Julia} implementation\footnote{The code is available at \url{https://github.com/kjellqvist/MinimaxAdaptiveControl.jl}.} of the proposed methods and design a controller that simultaneously stabilizes $G_\text{mp}$ and $G_\text{nmp}$,
	\begin{subequations}\label{eq:intro:example}
	\begin{align}
		G_\text{mp}(z) & = \frac{z_0z - 1/z_0}{(z - 1)^2},\label{eq:intro:mp} \\
		G_\text{nmp}(z) & =  \frac{z/z_0 - z_0}{(z - 1)^2}.\label{eq:intro:nmp}
	\end{align}
\end{subequations}
	Here $z$ is the complex frequency, and $z_0$ is $1.01$.
	Both systems are unstable, with a double pole in $1/(z-1)$.
	$G_\text{mp}$ has a minimum phase zero at $1/z_0^2$ and $G_\text{nmp}$ has a nonminimum-phase zero at $z_0^2$.
	The results are presented in Section~\ref{sec:examples:output}.
\subsection{Related work}\label{sec:related}
\subsubsection{Minimax Control}
Witsenhausen introduced minimax control in his thesis~\cite{Witsenhausen1966Minimax} as a decision-theoretic approach to control of uncertain dynamical systems.
Bertsekas and Rhodes~\cite{Bertsekas1973Sufficient} showed that the optimal controller can be decomposed into an estimator and an actuator. The optimal estimator can be expressed as a function of the observations, a so-called ``sufficiently informative function''.
$\Hinf$-control, both in the linear-quadratic case~\cite{Basar95Hinf} and in the nonlinear case~\cite{James95Robust}, has been formulated in terms of minimax control.
Recently, Goel and Hassibi \cite{Goel2023Competitive} showed that minimizing the regret or competitive ratio compared to an acausal controller with access to the future disturbance trajectory may provide excellent nominal performance at only a small robustness expense and that optimal controller synthesis can be reformulated as standard $\Hinf$ controller synthesis.
Karapetyan et. al~\cite{Karapetyan2022OnRegret} studied the suboptimality of $\Hinf$ control in the same setting.

The term ``minimax adaptive'' was introduced in~\cite{Basar1994MinimaxAdaptive} as a term for controllers that minimize the worst-case performance realization for systems with parametric uncertainty.
The authors considered continuous-time state-feedback control of systems with uncertain but constant parameters and showed that the cost can be rewritten in terms of the least-squares estimate of the parameters.
The reformulation using the least-squares estimate of the parameters is viable in the state-feedback case because it is sufficient to reconstruct the worst-feasible realization consistent with a model hypothesis and data, which James and Baras~\cite{James95Robust} showed is an information state---or informative statistic in Striebel's terms~\cite{Striebel1965Sufficient}.
This information state has a recursive formulation and is generally not finite-dimensional.
This statistic corresponds to our function $r$ in \eqref{eq:known_dynamics:r} and has a finite-dimensional representation.
Pan and Ba\c{s}ar generalized the results to of nonlinear SISO systems on ``parametric strict-feedback form'' in~\cite{Pan1998Adaptive}.

Vinnicombe~\cite{Vinnicombe2004} studied scalar systems where the parameters' signs are unknown and provided an explicit suboptimal controller based on certainty equivalence control with the least-squares parameter estimate.
Megretski and Rantzer~\cite{Megretski2003Bounds} provides lower bounds on the achievable $\ell_2$-gain for scalar systems with an uncertain pole belonging to an interval.
Rantzer extended Vinnicombe's result to higher-order systems where the state matrix has an unknown sign~\cite{Rantzer2020Minimax} and to finite sets of linear systems assuming full state measurements~\cite{Rantzer2021L4DC}.
A sufficient condition for finite $\ell_2$-gain is formulated as bilinear matrix inequalities, and a controller is obtained from the solution.

Cederberg \textit{et al.}~\cite{Cederberg2022Synthesis} proposed linearizing Rantzer's inequalities to improve performance iteratively, Bencherki and Rantzer~\cite{Bencherki2023Robust} gave conditions under which a solution to the inequalities is guaranteed to exist and Renganathan \textit{et al.}~\cite{Renganathan2023Online} studied empirical performance.
Kjellqvist generalized the framework to nonlinear systems where the preimage of the output under the measurement function is a finite set~\cite{Kjellqvist2024Minimax} assuming noise-free measurements.
Kjellqvist and Rantzer~\cite{Kjellqvist2022Noisy} previously extended Vinnicombe's~\cite{Vinnicombe2004} controller to the one-dimensional output-feedback case.

\subsubsection{Dual Control}
The controllers in Section~\ref{sec:explicit} are dual controllers.
In the nomenclature of Filatov and Unbehauen's survey\cite{Filatov2000Survey}, they are \emph{implicit} dual controllers, as they are suboptimal solutions to dual control problems.
Duality here is in the sense of Feldbaum's observation: that optimal controllers for uncertain nonlinear systems tend to have both regulating and experimenting mechanisms~\cite{Feldbaum1963Dual}.
This duality is known as the \emph{exploration-exploitation trade-off} in the reinforcement learning literature~\cite{Sutton1998Reinforcement}.
For further reading on dual control, see the surveys by Wittenmark~\cite{Wittenmark1995Dual}, Filatov and Unbehauen~\cite{Filatov2000Survey} and Mesbah~\cite{Mesbah2018Stochastic}.
\subsubsection{Supervisory Control \& Multiple-Model Adaptive Control}
Supervisory control, or multiple-model adaptive control, is a controller architecture where a supervisor selects a controller from a set of candidate controllers~\cite{Hespanha2001Tutorial}. 
Supervisory controllers typically come in two flavors: 
\emph{Estimator-based}, where each model has an associated estimator and control law, and the supervisor selects the model based on the estimator's output~\cite{Buchstaller2016Framework}, and \emph{controller-based}, where the supervisor selects among control laws by disqualifying controllers that violates assumed performance guarantees~\cite{Safonov1997Unfalsified, Patil2022Unfalsified}.
Our certainty-equivalence controllers can be seen as an instance of the former.

Switching, even among stable subsystems, may induce instability\cite{Liberzon2003Switching}, so the supervisor must ensure that the switching is safe.
In estimator-based frameworks, this typically translates to dwell-time constraints.
In the controller-based framework, this can be achieved by hysteresis in switching out underperforming controllers~\cite{Battistelli2010Stability}.
This switching restriction relates to the separation of time scales exploited in Ljung's averaging arguments~\cite{Ljung1977Analysis}, and to the difficulties of fast adaptation~\cite{Anderson2005Failures}.
The periodicity in our certainty-equivalence controllers can be seen as a form of dwell-time constraint.

\subsubsection{Learning}
Recently, there has been a surge of interest in the intersection of learning theory and control.
Advances in, for example, high-dimensional statistics~\cite{Tsiamis2023Statistical} and online convex optimization~\cite{Hazan2016Online} have provided tools for the design of new adaptive algorithms and analysis of achievable performance.
Most work has focused on relating the asymptotic scaling of performance bounds to assumptions on the size and number of uncertain parameters. It is based on assumptions about the statistics of the exogenous signals.
We now provide a brief overview of the works most closely related to our own.

Agarwal \textit{et al.}~\cite{Agarwal2019Online} considered control of linear systems with known dynamics, where the control objective and disturbances were adversarially chosen, as in our Problem~\ref{prob:known_dynamics}.
In contrast to our work, the cost functions were time-varying, revealed sequentially after actuation, and did not include a disturbance term.
Instead of a soft constraint, they assumed that the disturbance was bounded.

Simchowitz \textit{et al.}~\cite{Simchowitz2020Nonstochastic} extended the results to output feedback and uncertain dynamics but relied on apriori knowledge of a stabilizing static output feedback controller and evaluations of the objective function. 

Ghai \textit{et al.}~\cite{Ghai2022Robust} considered online control with model misspecification assuming perfect state measurements and provided an adaptive controller with bounded $\ell_2$-gain. 
The bound is asymptotic and scales with the number of uncertain parameters and the size of the model mismatch.
In contrast our $\ell_2$-bound is specific to the problem instance, and the user specifies precisely what parameters are uncertain.
	Lee \textit{et al.}~\cite{Lee2024Nonasymptotic} studied the regret of certainty-equivalence controllers with normally distributed exploration for approximately parameterized linear systems.
The approximation allows the user to inject prior knowledge and identify a reduced number of parameters, but means that the true system lies outside the space spanned by the parameters.
They showed an improvement over black-box adaptation in the small-data regime, as long as the model misspecification and the number of parameters are few.

\subsection{Notation}\label{sec:notation}
The set of $n\times m$ matrices with real coefficients is denoted $\R^{n \times m}$. 
The transpose of a matrix $A$ is denoted $A^\tran$. 
The space of real symmetric matrices in $\R^{n \times n}$ is denoted $\S^n$.
For a symmetrix matrix $H \in \S^{n + m}$ with blocks
\[
	H = \bmat{
		H^{11} & H^{12} \\
		H^{21} & H^{22}
	},
\]
we denote the Schur complements of $H^{11}$ and $H^{22}$ in $H$ by 
\[
	\begin{aligned}
		H / H^{22}  & = H^{11} - H^{12}H^{-22}H^{21},\\
		H / H^{11}  & = H^{22} - H^{21}H^{-11}H^{12}.
\end{aligned}
\]
$H^{-ii}$ denotes the inverse of $H^{ii}$, if it exists.
For a symmetric matrix $H \in \S^{n_1 + \cdots n_m}$ with blocks
\[
	H = \bmat{
		H^{11} & \cdots & H^{1m} \\
		\vdots & \ddots & \vdots \\
		H^{m1} & \cdots & H^{mm}
	},
\]
and vectors $x_1 \in \R^{n_1}, \ldots, x_m \in \R^{n_m}$, we define the quadratic form
\[
	\quadform_H(x_1, \ldots, x_m) = \sum_{i, j} x_i^\tran H^{ij} x_j.
\]
We write $H \succeq 0$ to indicate that $H$ is positive semidefinite and $H \succ 0$ to indicate that $H$ is positive definite.
If $H \succ (\succeq) 0$, we sometimes write $|x|^2_H$ to emphasize that $ |x|^2_H = \quadform_H(x)$ is a (semi) norm.
The standard euclidean norm is denoted $|x| = \sqrt{\sigma_I(x)}$, and by extension $|(x_1, \ldots, x_m)|^2 = \sum_{i=1}^m |x|^2$.
We refer to the value of a signal $w$ at time $t$ as $w_t$ and use the shorthand notation $\seq{w}{0}{t}$ for the sequence $(w_0, w_1, \ldots, w_t)$.
We sometimes use asterisks in matrix expressions to denote elements implied by symmetry.
For two sets, $X$ and $Y$, we denote the set of functions from $X$ to $Y$ by $Y^X$.

The \emph{convergence} and \emph{boundedness} of any sequence of functions $f_n : X \to Y$ in this paper are interpreted pointwise.

\section{Exact Analysis}\label{sec:exact}
\subsection{Principal Problem}\label{sec:principal}
This section introduces the principal problem of this paper and presents theory on minimax dynamic programming and value iteration.
In Sections~\ref{sec:state_feedback} and \ref{sec:output_feedback}, we show how to reduce state feedback and output feedback adaptive control to the principal problem.

\begin{problem}[Principal problem]\label{prob:known_dynamics}
	Let $\hat A \in \R^{n_z \times n_z}$, $\hat B \in \R^{n_z \times n_u}$, $\hat G \in \R^{n_z \times n_d}$, $z_0 \in \R^{n_z}$ and let $\mathcal M \subset \S^{n_z + n_u + n_d}$ be a compact set whose members, $H$, satisfy
	\begin{equation}\label{eq:known_dynamics:Q_constraints}
		\begin{aligned}
			H & = \bmat{
			H^{zz} & H^{zu} & H^{zd} \\
			H^{uz} & H^{uu} & H^{ud} \\
			H^{dz} & H^{du} & H^{dd}
		},
		& \begin{aligned}
			H^{dd} & \prec 0, \\
			H / H^{dd} & \succeq 0,\\
			[H / H^{dd}]^{uu} & \succ 0.
		\end{aligned}
	\end{aligned}
	\end{equation}
	Compute
	\begin{equation}\label{eq:known_dynamics:cost}
	\inf_\mu\sup_{H \in \modelset, d, N} \underbrace{\sum_{t=0}^{N-1} \quadform_H(z_t, u_t, d_t)}_{\hat J_\mu^N(z_0, H, d)},
	\end{equation}
	where $N \geq 0$ and the sequences, $\seq{z}{0}{N},\ \seq{u}{0}{N-1}$ are generated by
	\begin{align}
		z_{t + 1} & = \hat A z_t + \hat B u_t + \hat G d_t, \quad t \geq 0 \label{eq:known_dynamics:dynamics}\\
		u_t & = \mu_t(\seq{z}{0}{t}, \seq{u}{0}{t-1},\seq{d}{0}{t-1}). \label{eq:known_dynamics:control_law}
	\end{align}
\end{problem}

Problem~\ref{prob:known_dynamics} concerns the upper value of a two-player zero-sum game, where the minimization is over the controller, $\mu$, and the maximization is over the disturbance, $d$, and the realization of the cost function, $H \in \modelset$.
If not for the uncertainty in the cost function, the problem would be a (nonstandard) linear-quadratic control problem, which is a well understood problem class~\cite{Basar95Hinf}.

The relation to adaptive control is as follows.
In state-feedback adaptive control, with dynamics of the form $x_{t+1} = Ax_t + Bu_t + w_t$, where $x_t$ is the state and $w_t$ is the disturbance and the pair $(A, B)$ is unknown, and quadratic stage costs, we let $d_t = x_{t+1}$ and $z_t = x_t$.
Substituting $w_t = d_t - Az_t - Bu_t$ and $z_t$ into the cost function gives dynamics of the form~\eqref{eq:known_dynamics:dynamics} and cost functions of the form~\eqref{eq:known_dynamics:cost}.
This is explained in more detail in Section~\ref{sec:state_feedback}.

For output-feedback adaptive control with a finite set of feasible models and quadratic stage costs, we quantify the worst-case accrued cost using one observer for each model.
The $z_t$ of Problem~\ref{prob:known_dynamics} is constructed by stacking the observer states, $d_t$ is the measured output, and $u_t$ is the control input.
The matrices $\hat A$, $\hat B$ and $\hat G$ corresponds to aggregating the observer dynamics.
We get one Hessian, $H$, for each model, expressing the past performance of the observer.
The reformulation of output-feedback adaptive control as an instance of Problem~\ref{prob:known_dynamics} is explained in Section~\ref{sec:output_feedback}.
The rest of this section is devoted to dynamic programming and value iteration.
\subsubsection{Dynamic Programming}
Define the functions $r_t : \modelset \to \R$ for $t = 0, 1, 2, \ldots$ by
\begin{equation}\label{eq:known_dynamics:r}
	r_{t}(H) = \sum_{s= 0}^{t-1}\quadform_H(z_s, u_s, d_s).
\end{equation}
Then $r_t$ satisfies the recursion
\begin{equation}
	r_{t+1}(H) = \underbrace{r_{t}( H) + \quadform_H(z_t, u_t, d_t)}_{f(r_t, z_t, u_t, d_t)(H)}.
\end{equation}
Although the controller does not know the realization of $H$, the functions $r_{t}$ are constructed of known quantities and can be computed by the controller at time t.

\begin{remark}
	The functions $r_t$ take the form $r_t(H) = \langle Z_t, H \rangle$.
	The positive semidefinite matrix $Z_t$ can be computed recursively by
	\[
		Z_{t + 1} = Z_t + \bmat{z_t\\ u_t \\ d_t} \bmat{z_t \\ u_t \\ d_t}^\tran, \quad Z_0 = 0.
	\]
	This means that the matrix $Z_t$ compresses the information of the past states, inputs, and disturbances into a single matrix.
	If the cardinality of $\modelset$ is large compared to the state dimension, $n_z$, then the matrix $Z_t$ can be used to reduce the computational complexity of the problem.
	If the model set is finite, then one can store the function values $r_t(H)$ for each $H \in \modelset$ in an array.
\end{remark}
For a function $V : \R^{n_z} \times \R^\modelset \to \R$, define the Bellman operators
\begin{subequations}\label{eq:known_dynamics:Bellman_both}
\begin{align}
	\B_u V(z, r) & = \max_d V\left(\hat Az + \hat Bu + \hat Gd, f(r, z, u, d)\right)\label{eq:known_dynamics:Bellman_u}, \\
	\B V(z, r) & = \min_u\B_u V(z, r), \label{eq:known_dynamics:Bellman}
\end{align}
\end{subequations}
and the value iteration
\begin{subequations}\label{eq:known_dynamics:value_iteration}
\begin{align}
	V_0(z, r) & = \max_{H\in \modelset}r(H),\\
	V_{k+1}(z, r) & = \B V_k(z, r).
\end{align}
\end{subequations}
We will consider control policies $\eta_t : \R^{n_z} \times \R^\modelset \to \R^{n_u}$ of the form
\begin{equation}\label{eq:known_dynamics:policy}
	u_t = \eta_t(z_t, r_t),
\end{equation}
and note that this policy is admissible as $r_t$ depends causally on the states, inputs, and measured disturbances.

\newcounter{prop:known_dynamics:summary}%
\begin{theorem}\label{thm:known_dynamics}
	The following facts holds for Problem~\ref{prob:known_dynamics}, the Bellman operator $\B$ in \eqref{eq:known_dynamics:Bellman} and the value iteration defined in \eqref{eq:known_dynamics:value_iteration}.
	\begin{enumerate}
		\item $\B$ is monotone: $V' \geq V \implies \B V' \geq \B V$.
		\item The value iteration is nondecreasing:  $V_{k + 1} \geq V_k$.
		\item The value iteration converges if, and only if, it is bounded.
		\item The value~\eqref{eq:known_dynamics:cost} is bounded for all $z_0 \in \R^{n_z}$ if, and only if, the value iteration converges.
			\setcounter{prop:known_dynamics:summary}{\value{enumi}}
	\end{enumerate}

If the value iteration converges to a limit $V_\star$, then
	\begin{enumerate}
		\setcounter{enumi}{\value{prop:known_dynamics:summary}}
		\item The value \eqref{eq:known_dynamics:cost} is equal to $V_\star(z_0, 0)$.
		\item $V_\star$ is a fixed point of $\B$, not necessarily unique.
		\item $V_\star$ is the minimal fixed point of $\B$ greater than $V_0$.
		\item The control law $\eta_\star$ defined as the minimizer in~\eqref{eq:known_dynamics:Bellman} achieves \\ $\B_{\eta_\star(z, r)}V_\star (z, r) = \B V(z, r)$, and the policy
	\[
		\mu_t(\seq{z}{0}{t},\seq{u}{0}{t-1}, \seq{d}{0}{t-1}) = \eta_\star(z_t, r_t)
	\]
	is optimal for Problem~\ref{prob:known_dynamics}.

	\end{enumerate}
\end{theorem}
\begin{proof}
	The proofs of the statements in the theorem are standard but we include them here for completeness.
	\begin{enumerate}
		\item As the maximization in \eqref{eq:known_dynamics:Bellman_u} and minimization in \eqref{eq:known_dynamics:Bellman} are monotone operations, so is their composition $\B$.
		\item Assume that $V_k \geq V_{k-1}$. By monotonicity of $\B$, $V_{k+1} = \B V_k \geq \B V_{k-1} = V_k$. We now consider the base case, and prove that $V_1 \geq V_0$. By the minmax inequality
			\begin{align*}
				V_1(z, r) & = \min_u\max_d V_0(\hat Az + \hat Bu + \hat Gd, f(r, z, u, d)) \\
					  & \geq \max_{H \in \modelset} \min_u \max_d \Big\{\sigma_H(\hat Az + \hat Bu + \hat Gd, u, d) \\
					  & \quad+ r(H)\Big\}
			\end{align*}
			By \eqref{eq:known_dynamics:Q_constraints}, $\max_d \sigma_H(z, u, d) \geq 0$ for all $H \in \modelset$, we have that $V_1(z, r) \geq \max_{H\in \modelset} r(H) = V_0(z, r)$.
		\item Pointwise convergence of the value iteration is equivalent to convergence of the monotone sequence of real numbers $V_0(z, r), V_1(z, r)$, \ldots.
	\end{enumerate}

	We first show that $V_k(z, r)$ converges if, and only if $V_k(z, 0)$ converges. Let $a \in \R$ be a constant. By induction $V_k(z, r + a) = V_k(z, r) + a$. Further, as $f$ and $V_0$ are monotone in $r$, so is $V_k$: $r' \geq r \implies V_k(z, r') \geq V_k(z, r)$.
			Thus, for all $k = 1, 2 \ldots$
			\[
				V_k(z, 0) + \min \{r\} \leq V_k(z, r) \leq V_k(z, 0) + \max \{r\}.
			\]
			For any $N$, standard dynamic programming arguments show that $V_N(z_0, 0) = \inf_\mu \sup_{d, H \in \modelset} \hat J_\mu^N(z_0, H, d)$, which is a lower bound for \eqref{eq:known_dynamics:cost}. Thus \eqref{eq:known_dynamics:cost} is bounded \emph{only if} the value iteration is bounded.

			Now, fix an $\epsilon > 0$ and let $\tilde \eta_t$ be a policy that achieves $\B_{\eta_t} V_\star(z, r) \leq V_\star(z, r) + \epsilon^{t + 1}$. 
			Then \eqref{eq:known_dynamics:cost} is bounded from above by the expression
			\begin{align*}
				\sup_{N, d, H \in \modelset} \hat J_{\tilde \eta}^N(z_0, H, d) & \leq \sup_N\left\{ V_\star(z, r) + \epsilon\sum_{t = 0}^N \epsilon^t \right\}\\
				= V_\star (z, 0) + \frac{\epsilon}{1 - \epsilon}.
			\end{align*}
			Thus we conclude statements 4, 5, and 8.

			As $V_\star$ is the limit of the value iteration, it is a fixed point of $\B$ (otherwise, the limit would not exist).
			To see that it is minimal, assume that $V'$ is another fixed point of $\B$ such that $V' \geq V_0$ but that $V'(z_0, r_0) < V_\star(z, r)$ for some $z_0$ and $r_0$.
			Define the function $\hat V_0$ by
			\[
				\hat V_0(z, r) = \begin{cases}
					V'(z, r) & \text{if } (z, r) = (z_0, r_0), \\
					V_0(z, r) & \text{otherwise}.
				\end{cases}
			\]
			Then $V_0 \leq \hat V_0 \leq V'$, so the value iteration $\hat V, \B \hat V, \ldots$ converges to some limit $\hat V_\star$. But by monotonicity, $V_\star \leq \hat V_\star \leq V'$, which is a contradiction.
\end{proof}

\subsection{State Feedback}\label{sec:state_feedback}
\begin{problem}[State-feedback minimax adaptive control]\label{prob:state_feedback}
	Let $Q \in \S^{n_x}$, $R \in \S^{n_u}$, be positive definite. 
	Given a compact set $\modelset \subset \R^{n_x \times n_x} \times \R^{n_x \times n_u}$, initial state $x_0 \in \R^{n_x}$, and a positive quantity $\gamma > 0$, compute
	\begin{equation}\label{eq:adaptive_state_feedback:cost}
		\inf_\mu\sup_{w, N, M \in \modelset} \underbrace{\sum_{t=0}^{N-1} \left(|x_t|^2_Q + |u_t|^2_R - \gamma^2|w_t|^2 \right)}_{J_\mu^N(x_0, M, w)}
	\end{equation}
	where $M = (A, B) \in \modelset$, $w_t \in \R^{n_w}$, $N \geq 0$, and the sequences $\seq{x}{0}{N},\ \seq{u}{0}{N-1}$ are generated by
	\begin{align}
		x_{t + 1} & = A x_t + B u_t + w_t, \quad t \geq 0 \label{eq:state_feedback:dynamics}\\
		u_t & = \mu_t(\seq{x}{0}{t}, \seq{u}{0}{t-1}). \label{eq:state_feedback:control_law}
	\end{align}
\end{problem}

The state-feedback minimax adaptive control problem, Problem \ref{prob:state_feedback}, is similar to a standard $\Hinf$ control problem, but differs in that the dynamics are uncertain and chosen by the adversary. 
The problem differs from the principal problem in that the realization of the objective function is known, but the dynamics are not.
To relate the two problems, we introduce $z_t = x_t$, $d_t = x_{t+1}$.
Substituting $w_t = d_t - Ax_t - Bu_t$ into the dynamics \eqref{eq:state_feedback:dynamics}, we get
\begin{equation}
	z_{t + 1} = \hat A z_t + \hat B u_t + \hat Gd_t, \\
\end{equation}
where $\hat A = 0$, $\hat B = 0$ and $\hat G = I$.
For $M = (A, B) \in \modelset$, let
\begin{equation}\label{eq:state_feedback:H}
	H_M = \bmat{
		Q & 0 & 0\\
		0 & R & 0\\
		0 & 0 & 0
	}
	- \gamma^2\bmat{
		-A^\tran \\ -B^\tran \\ I
	}
	\bmat{
		-A^\tran \\ -B^\tran \\ I
	}^\tran.
\end{equation}
Then, the objective \eqref{eq:adaptive_state_feedback:cost} becomes
\begin{equation}\label{eq:state_feedback:Riccati_operator}
	\inf_\mu\sup_{d, N, M \in \modelset} \sum_{t=0}^{N-1} \quadform_{H_M}(z, u, d).
\end{equation}

\[
	H_M / H_M^{dd} = \bmat{Q & 0 \\ 0 & R} \succeq 0.
\]
Finally, note that $H_M$ fulfills \eqref{eq:known_dynamics:Q_constraints} as $R \succ 0$.
We summarize the result in the following theorem.
\begin{theorem}[State-feedback reduction]
	The value of Problem \ref{prob:state_feedback} is equal to the value of Problem \ref{prob:known_dynamics} with $\modelset = \{H_M : M \in \modelset\}$ with $H_M$ as in \eqref{eq:state_feedback:H}, $z_t = x_t$, $d_t = x_{t+1}$, and $\hat A = 0$, $\hat B = 0$, and $\hat G = I$.
Further, given a policy $\hat \mu : (\seq{z}{0}{t},\seq{u}{0}{t-1}, \seq{d}{0}{t-1}) \mapsto u_t$, 
	\[
		\sup_{d, N, H \in \widehat \modelset}\hspace{-1em} \hat J_{\hat \mu}^N(z_0, H, d) =\hspace{-1em} \sup_{w, N, M\in \modelset}\hspace{-1em} J^N_{\mu}(x_0, M, w),
	\]%
	where%
	\begin{equation*}
		\mu(\seq{x}{0}{t}, \seq{u}{0}{t}) \\
		= \hat \mu(\seq{z}{0}{t}, \seq{u}{0}{t-1}, \seq{d}{0}{t-1})
	\end{equation*}
	is feasible for Problem~\ref{prob:output_feedback}.
\end{theorem}

\subsection{Output Feedback}\label{sec:output_feedback}
This section presents how to rewrite the output-feedback minimax adaptive control problem formalized below as an instance of the principal problem, Problem \ref{prob:known_dynamics}.
\begin{problem}[Output-feedback minimax adaptive control]\label{prob:output_feedback}
	Let $Q \in \S^{n_x}$, $R \in \S^{n_u}$, be positive definite. 
	Given a compact set $\modelset \subset \R^{n_x \times n_x} \times \R^{n_x \times n_u} \times \R^{n_x \times n_w} \times \R^{n_y \times n_x} \times \R^{n_y \times n_v}$, and a positive quantity $\gamma > 0$ , consider 
	\begin{equation}\label{eq:output_feedback:objective_function}
		\underbrace{- |x_0 - \hat x_ 0|^2_{S_{M, 0}} 
	+ \sum_{t=0}^{N-1} \left(|x_t|^2_Q + |u_t|^2_R - \gamma^2|(w_t, v_t)|^2 \right)}_{J_\mu^N(\hat x_0, M, y, w, v, x_0)}.
\end{equation}
	where $S_{M,0} \succeq 0$, $M = (A, B, G, C, D) \in \modelset$, $w_t \in \R^{n_w}$, $v_t \in \R^{n_v}$, $N \geq 0$, and the sequences, $\seq{x}{0}{N}$, $\seq{y}{0}{N-1}$ and $\seq{u}{0}{N-1}$ are generated by%
\begin{subequations}\label{eq:output_feedback:system}
	\begin{align}
		x_{t + 1} & = A x_t + B u_t + Gw_t, \quad t \geq 0 \label{eq:output_feedback:dynamics}\\
		y_t & = C x_t + D v_t, \label{eq:output_feedback:output}\\
		u_t & = \mu_t(\seq{y}{0}{t-1}, \seq{u}{0}{t-1}). \label{eq:output_feedback:control_law}
	\end{align}
\end{subequations}

	Compute
	\begin{equation}\label{eq:output_feedback:objective}
		\inf_\mu\sup_{y, M \in \modelset, N, w, v}J^N_\mu(\hat x_0, M, y, w, v).
	\end{equation}
\end{problem}

\begin{remark}
	We assume that all members of $\modelset$ have the same order, i.e., $n_x$, $n_w$ and $n_v$ are constant for all $M \in \modelset$. This is for notational simplicity only and is not necessary for the theoretical development in this section.
\end{remark}
We make the following assumptions on the problem parameters.
\begin{assumption}[Problem parameters]\label{assumption:output_feedback:problem}
	For each $M \in \mathcal M$,
	$\bmat{
			A & G \\ 
		}
	$
	and $D$ have full column rank, and $\bmat{A^\tran & \sqrt{Q}}^\tran$ has full row rank.
\end{assumption}
We follow the approaches of~\cite{James95Robust} and \cite{Basar95Hinf}, and split the optimization problem \eqref{eq:output_feedback:objective} into three steps,
\begin{multline}\label{eq:output_feedback:W}
	\inf_\mu \sup_{y, M\in\modelset, w, v, N, x_0} J_\mu^N(\hat x_0, M, y, w, v, x_0) =\\
	\inf_\mu \sup_{y, M \in \modelset, N}\hspace{-.4em}\sup_{x} \underbrace{\left[\sup_{w, v, x_0} \{ J_\mu^N(\hat x_0, M, y, w, v, x_0) : x_N = x\}\right]}_{W_M^N(x, \seq{u}{0}{N-1}, \seq{y}{0}{N-1})}.
\end{multline}
The supremum in $W_M^N(x, \seq{u}{0}{N-1}, \seq{y}{0}{N-1})$ is taken subject to the inputs $\seq{w}{0}{N-1}$, $\seq{v}{0}{N-1}$ and $\seq{u}{0}{N-1}$, the observed output $\seq{y}{0}{N-1}$, final state $x_N = x$ and model $M$ being feasible.
Feasibility means that the state sequence $\seq{x}{0}{N}$ is generated by the dynamics \eqref{eq:output_feedback:dynamics} and the output sequence $\seq{y}{0}{N-1}$ is generated by the output equation \eqref{eq:output_feedback:output} under the model $M$.

Note that $W^N_M(x, u, y)$ is a function of the trajectory $\seq{u}{0}{N-1}$ and not the control law $\mu$.
This is because the outer optimization steps determine the control law and the output sequence.
The control law and the output sequence in turn determine the control signal trajectory.

The reformulation has two major benefits. 
The first is that minimizing over $\mu_t$ and maximizing over $\seq{y}{0}{t-1}$ commutes as $\mu_t$ is a function of $\seq{y}{0}{t-1}$.
This interchange of extremization leads to a sequential optimization problem that can be solved by backwards dynamic programming.
The second benefit is that $W_M^N(x, \seq{u}{0}{N-1}, \seq{y}{0}{N-1})$ can be characterized using standard forward Riccati recursions and observer equations of $\Hinf$-control theory.

The rest of the section is organized as follows.
Section~\ref{sec:output_feedback:known_dynamics} shows how to rewrite Problem \ref{prob:output_feedback} as an instance of Problem \ref{prob:known_dynamics} in the case where the dynamics are known.
Section~\ref{sec:output_feedback:unknown_dynamics} modifies the approach to the case where $\modelset$ is a finite set of models.

\subsubsection{Known Dynamics}\label{sec:output_feedback:known_dynamics}
In this case, $\modelset = \{M\}$ for some $M = (A, B, G, C, D)$, and we will drop the subscript $M$ from the notation.
The value \eqref{eq:output_feedback:objective} is then equal to
\[
	\inf_\mu \sup_{N, \seq{y}{0}{N-1}, x} W^N(x, \seq{u}{0}{N-1}, \seq{y}{0}{N-1}).
\]

Consider the forward Riccati recursions:
\begin{equation}\label{eq:certain:output_feedback:forward_riccati}
\begin{aligned}
	S_{t+1} & = (A X_t^{-1} A^\tran + \gamma^{-2}GG^\tran)^{-1}, \\
	X_t & = S_t + \gamma^2 C^\tran(DD^\tran)^{-1} C - Q, \\
	L_t & = \gamma^2 A X_t^{-1} C^\tran(DD^\tran)^{-1},
\end{aligned}
\end{equation}
where $S_t \in \S^{n_x}$, $X_t \in \S^{n_y}$, $L_t \in \R^{n_x \times n_y}$.
The $\Hinf$-observer states obey the dynamics
\begin{equation}\label{eq:Hinf_observer}
	\hat x_{t + 1} = A \hat x_{t} + Bu_t + L_t(y_t - C \hat x_{t}) +AX_t^{-1}Q\hat x_t.
\end{equation}
The initial $S_0$ and $\hat x_0$ are provided by the designer in~\eqref{eq:output_feedback:objective_function}.
The following lemma summarizes the recursive computation of $W^N$, we refer the reader to~\cite[Chapter 6]{Basar95Hinf} for a proof.
\begin{lemma}\label{lemma:prel:output_feedback:W}
	For Problem~\ref{prob:output_feedback} with $\modelset = \{(A, B, C, G, D)\}$, let $S_{t + 1} \in \S^{n_x}$, $X_t \in \S^{n_y}$, and $L_t \in \R^{n_x \times n_y}$ be defined recursively for $t = 0, 1, \ldots$ by \eqref{eq:certain:output_feedback:forward_riccati}.
	If $S_t\succ Q$ for all $t = 0, \ldots, N$, then $W^N$ defined in \eqref{eq:output_feedback:W} satisfies
	\begin{equation}\label{eq:certain:output_feedback:W^N}
	W^N(x, \seq{u}{0}{N-1}, \seq{y}{0}{N-1}) = -|x - \hat x_N|^2_{S_N} + r_N.
	\end{equation}
	The observer states $\hat x_t$ follow the observer equation~\eqref{eq:Hinf_observer}, $r_0 = 0$ and
	\begin{multline}
		r_{t + 1}  = r_t + |\hat x_{ t}|^2_{(Q^{-1} - S_t^{-1})^{-1}} + |u_t|^2_R \\
			 - |y_t - C (S - Q)^{-1}S\hat x_{ t}|^2_{(DD^\tran / \gamma^2 + C(S - Q)^{-1}C^\tran)^{-1}}.
	\end{multline}
If, for any $t \geq 0$, $S_t \nsucc Q$, then the value \eqref{eq:output_feedback:objective} is unbounded.
\end{lemma}

It is not obvious how the designer should choose the initial $S_0$.
The Riccati recursions~\eqref{eq:certain:output_feedback:forward_riccati} are known to admit positive definite fixed points $S$ if $\gamma$ is sufficiently large, and the minimal fixed point leads to stable observer dynamics~\eqref{eq:Hinf_observer}.
Thus, we can choose $S_0 = S$, a stabilizing positive definite fixed point of the Riccati recursions, and will assume so for the rest of the section.
\begin{assumption}\label{assumption:output_feedback:fixed_point}
	$\gamma$ is large enough so that there exists a stabilizing positive definite fixed point of the Riccati recursions~\eqref{eq:certain:output_feedback:forward_riccati}.
	Denote by $S$ the minimal fixed point and assume $S_0 = S_M$.
\end{assumption}

\begin{remark}
Note that under Assumption~\ref{assumption:output_feedback:fixed_point}%
	\[
		r_{t} = \sum_{s = 0}^{t-1} \quadform_{\hat Q}(\hat x_s, u_s, y_s)
	\]%
	where
	\begin{equation}\label{eq:certain:output_feedback:Q_tau}
		\begin{aligned}
		 & 	\hat Q^{xx} =  SX^{-1} S - S \quad \hat Q^{xu}= 0, \\
		 &	\hat Q^{uu} = R, \quad
		 	\quad \hat Q^{yy} = -(DD^\tran / \gamma^2 + C(S- Q)^{-1}C^\tran)^{-1}, \\
		 &	\hat Q^{uy} = 0, \quad
		 \hat Q^{xy}= \gamma^2 SX^{-1}C^\tran(DD^\tran)^{-1}.
		\end{aligned}
	\end{equation}
\end{remark}

\begin{remark}
The expression of $r_{t+1}$ is unnecessary when one employs the certainty equivalence principle of~\cite{Basar95Hinf}, but is crucial to our theory of adaptive control.
\end{remark}

\begin{theorem}\label{thm:certain:output_feedback:reduction}
	For Problem~\ref{prob:output_feedback} with $\modelset = \{(A, B, C, G, D)\}$, under assumptions \ref{assumption:output_feedback:problem} and \ref{assumption:output_feedback:fixed_point}, let $\hat A = A X^{-1} S$, $\hat B = B$, $\hat G = L$, where $(S, X, L)$ is a fixed point of the Riccati recursions~\eqref{eq:certain:output_feedback:forward_riccati}.
Then, the optimal value of Problem~\ref{prob:output_feedback} is equal to the optimal value of Problem~\ref{prob:known_dynamics} where $H$ is replaced by $\hat Q$ in~\eqref{eq:certain:output_feedback:Q_tau} and $z_t = \hat x_t$.
Further, given a policy $\hat \mu : (\seq{z}{0}{t},\seq{u}{0}{t-1},\seq{y}{0}{t-1})) \mapsto u_t$, 
	\[
		\sup_{d, N} \hat J_{\hat \mu}^N(z_0, H, d) =\hspace{-.5em} \sup_{y, N, w, v, x_0}\hspace{-.5em} J^N_{\mu}(\hat x_0, M, y, w, v, x_0),
	\]
	where 
	\begin{equation*}
		\mu(\seq{y}{0}{t-1}, \seq{u}{0}{t-1})
		= \hat \mu(\seq{z}{0}{t}, \seq{u}{0}{t-1}, \seq{y}{0}{t-1})
	\end{equation*}
	is feasible for Problem~\ref{prob:output_feedback}.
\end{theorem}

\begin{proof}
		From Lemma~\ref{lemma:prel:output_feedback:W}, we know that for a fixed policy $\mu$ and horizon $N \geq 0$, the value of the inner optimization problem in~\eqref{eq:output_feedback:objective},
\begin{multline*}
	\sup_{w, v} \left\{\sum_{t = 0}^{N-1}\left(|x_t|^2_Q + |u_t|^2_R - \gamma^2|(w_t, v_t)|^2\right) - |x_0 - \hat x_0|^2_{S_0} \right\} \\
	=
	\underbrace{\sup_{\seq{y}{0}{N-1}, x}\left\{
		-|x - \hat x_N|_{S_N} + \sum_{t = 0}^{N - 1}\quadform_{\hat Q}(\hat x_t, u_t, y_t)
\right\}}_{J^N_\mu(\hat x_0)},
\end{multline*}
	where $\hat Q$ is defined in~\eqref{eq:certain:output_feedback:Q_tau} and the sequences $\hat x_t$ and $u_t$ are generated by~\eqref{eq:Hinf_observer} and~\eqref{eq:output_feedback:control_law}, respectively.
	As $S_N$ is positive definite, the unique maximizing argument is $x_\star = \hat x_N$.
	By assumption, $S_0$ is a fixed point of the Riccati recursion~\eqref{eq:certain:output_feedback:forward_riccati} and therefore $S_t = S = S_0$.
	The corresponding matrices $\hat Q_t$, $X_t$, $L_t$, $\hat A$, $\hat B$ and $\hat G$ are also stationary.
	Finally, as $\hat x_t$ is a function of $\seq{y}{0}{t-1}$ and $\seq{u}{0}{t-1}$, the sets of feasible controllers $\mu$ in Problem~\ref{prob:output_feedback} and $\hat \mu$ in this theorem are equal.
	We conclude that the optimal values are equal.
\end{proof}
\subsubsection{Main Result: Output Feedback Adaptive Control}\label{sec:output_feedback:unknown_dynamics}
We now consider the case where the dynamics of the system are unknown and the controller must adapt to the system, but assume that the model set $\modelset$ is finite.
Let $\hat x_{M, t}$, $\hat A_M$, $\hat B_M$, $\hat G_M$ and $\hat Q_M$ be as in Theorem~\ref{thm:certain:output_feedback:reduction} for each $M \in \modelset$ and define
\begin{equation}\label{eq:output_feedback:reduction}
	z_t = \bmat{\hat x_{1, t}\\ \vdots \\ x_{|\modelset|, t}}, \quad
	\hat B = \bmat{\hat B_1 \\ \vdots \\ \hat B_{|\modelset|}}, \quad
	\hat G = \bmat{\hat G_1 \\ \vdots \\ \hat G_{|\modelset|}},
\end{equation}
\[
	\hat A = \bdiag\{A_MX_M^{-1} S_M: M \in \modelset\},
\]
and $d_t = y_t$. 
Further, let $H_i \in \S^{|\modelset|n_x \times n_u \times n_y}$ be given by
\vspace{2em}
{\small
\begin{equation}\label{eq:output_feedback:hat_Q}
	\qquad \left[{\color{gray}
		\begin{array}{ccccccc|c|c}
			\color{gray}0 & \cdots & 0 &\tikz[remember picture]\node[inner sep = 0pt] (coli) {\color{black}0}; &0 & \cdots  & 0 & 0 & 0 \\
			\vdots & \ddots &&{\color{black} \vdots}&&& \vdots & \vdots & \vdots \\
			0 & & 0 & {\color{black}0} & 0 & \cdots & 0 & 0 & 0 \\
		
		\tikz[remember picture]\node[inner sep = 0pt] (rowi) {\color{black}0}; &{\color{black}\cdots} & {\color{black}0} & {\color{black}\hat Q_i^{xx}} & {\color{black}0} & {\color{black}\cdots} & {\color{black}0} &{\color{black}\hat Q_i^{xu}}& {\color{black}\hat Q_i^{xy}}\\
		0 & & 0 & {\color{black}0} & 0 & \cdots & 0 & 0 & 0 \\
		\vdots & &&{\color{black} \vdots} & &  \ddots & \vdots & \vdots & \vdots \\
			0 & \cdots & 0 & {\color{black}0} & 0 & \cdots & 0 & 0 & 0 \\
			\hline
			0 & \cdots & 0 & {\color{black}\hat Q_i^{ux}} & 0 & \cdots & 0 & \hat Q_i^{uu} & \hat Q_i^{uy} \\
			\hline
			0 & \cdots & 0 & {\color{black}\hat Q_i^{yx}} & 0 & \cdots & 0 & \hat Q_i^{yu} & \hat Q_i^{yy}
		\end{array}
		}
	\right].
\end{equation}
}

\begin{tikzpicture}[remember picture, overlay]
	\node [left = 1em of rowi, align = center] {row $i$};	
	\node [above= 1.2em of coli, align = center] (colitext) {column $i$};	
	\draw[->] ([yshift=0.4em]coli.north) -- (colitext);
\end{tikzpicture}

The following theorem shows that the output-feedback adaptive control problem can be reduced to an instance of Problem~\ref{prob:known_dynamics}.
\begin{theorem}[Reduction]\label{thm:output_feedback:reduction}
	Under Assumptions~\ref{assumption:output_feedback:problem} and \ref{assumption:output_feedback:fixed_point}, the optimal value of Problem~\ref{prob:output_feedback} is equal to the optimal value of Problem~\ref{prob:known_dynamics} where $\modelset$ is replaced by $\{H_M : M \in \modelset\}$ with $H_M$ as in~\eqref{eq:output_feedback:hat_Q} and $\hat A$, $\hat B$ and $\hat G$ in~\eqref{eq:output_feedback:reduction}, and $z_0$ in~\eqref{eq:output_feedback:reduction}.
	Further, given a policy $\hat \mu : (\seq{z}{0}{t},\seq{u}{0}{t-1},\seq{y}{0}{t-1}) \mapsto u_t$, 
	\[
		\sup_{d, N, H \in \widehat \modelset} \hat J_{\hat \mu}^N(z_0, H, d) =\hspace{-1.5em} \sup_{y, M\in \modelset, N, w, v, x_0}\hspace{-1.5em} J^N_{\mu}(\hat x_0, M, y, w, v, x_0),%
	\]%
	where%
	\begin{equation*}
		\mu(\seq{y}{0}{t-1}, \seq{u}{0}{t-1})
		= \hat \mu(\seq{z}{0}{t}, \seq{u}{0}{t-1}, \seq{y}{0}{t-1})
	\end{equation*}
	is feasible for Problem~\ref{prob:output_feedback}.
\end{theorem}
\begin{proof}
	The proof is identical to that of Theorem~\ref{thm:certain:output_feedback:reduction}, but with the additional steps mentioned below.
	With $W$ as in \eqref{eq:output_feedback:W}, we have that the cost~\eqref{eq:output_feedback:objective} is equal to
	\[
		\inf_\mu \sup_{M, N, \seq{y}{0}{N-1}, x} W_M^N(x, \seq{u}{0}{N-1}, \seq{y}{0}{N-1}).
	\]
	Let $\hat S_M \in \R^{|\modelset|n_x \times |\modelset|n_x}$ be the matrix that is zero except on the $M$-th diagonal block, which is equal to $S_M$.
	By Lemma~\ref{lemma:prel:output_feedback:W}, we have that
	\begin{align*}
		& \sup_x W_M^N(x, \seq{u}{0}{N-1}, \seq{y}{0}{N-1}) =\sup_x\{ -|x - \hat x_{M, N}|^2_{\hat S_M} + r_{M, N}\} \\
		& = \sum_{t = 0}^{N-1}\quadform_{\hat Q_M}(\hat x_{M, t}, u_t, y_t) 
								 = \sum_{t = 0}^{N-1}\quadform_{H_M}(z_t, u_t, y_t).
	\end{align*}
	where $\hat Q_M$ is as in~\eqref{eq:certain:output_feedback:Q_tau}, $z_t$ is as in~\eqref{eq:output_feedback:reduction}, and $H_M$ is as in~\eqref{eq:output_feedback:hat_Q}.
\end{proof}

\begin{figure}
	\centering
	\includegraphics[width=\linewidth]{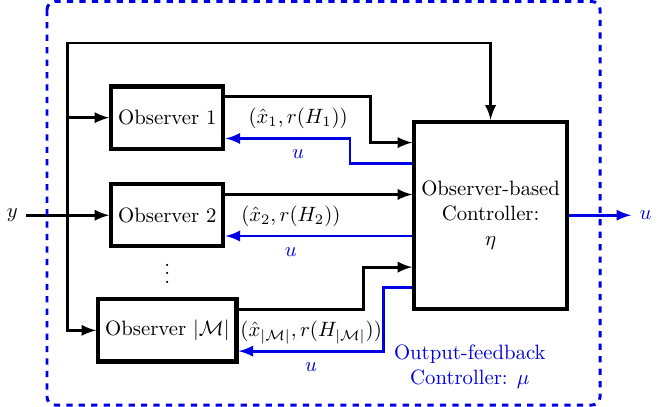}
	\caption{Observer-based adaptive controller resulting from combining Theorem~\ref{thm:output_feedback:reduction} and Theorem~\ref{thm:known_dynamics}.
	}
	\label{fig:obscontrol}
\end{figure}

Together with Theorem~\ref{thm:known_dynamics}, Theorem~\ref{thm:output_feedback:reduction} can be realized by a bank of observers as in Fig.~\ref{fig:obscontrol}.

\section{Explicit Controller Synthesis with Performance Bounds}\label{sec:explicit}
\subsection{Bellman Inequalities}\label{sec:inequalities}
Sometimes, it may be difficult to compute the recursion~\eqref{eq:known_dynamics:value_iteration} or to find the minimal nonnegative fixed point of the Bellman operator $\B$ in \eqref{eq:known_dynamics:Bellman} and the corresponding optimal control law.
This is typically the case in our setting. 
Some exceptions are the case of $\modelset$ being a singleton or the case where the uncertainty set $\modelset$ contains an element $H_M$ that dominates all other models.
By dominate, we mean that a control policy that is optimal for $H_M$ achieves a lower cost for all other models in $\modelset$.

This section presents theory related to bounding the value function that relies on periodic compensation.
As we will see in Section~\ref{sec:examples}, delays in the control signal significantly complicates the computation of the value-function approximation.
Periodic compensation is a powerful method to handle this problem.

Let $\tau$ be a positive integer. 
We model $\tau$-periodic compensation as a control law that contains a supervisor, $\bar \eta$, that periodically selects a sequence of $\tau$ control components,
\begin{equation}\label{eq:periodic_compensation}
	\bar \eta : (z_{\tau k}, r_{\tau k}) \mapsto (\bar \eta_{\tau k}, \bar \eta_{\tau k + 1}, \ldots, \bar \eta_{\tau k + \tau - 1}).
\end{equation}
During this period, the control signal is computed by the component control laws
\begin{equation}\label{eq:periodic_compensation:control}
		u_{\tau k + s} = \bar \eta_{\tau k + s}(z_{\tau k + s}, r_{\tau k + s}).
\end{equation}

The connection to dynamic programming lies in compositions of the Bellman operator.
The periodic versions of the operators $\B$ and $\B_u$ in \eqref{eq:known_dynamics:Bellman_both} are
\begin{subequations}\label{eq:periodic_compensation:operators}
	\begin{align}
		\B^\tau \bar V & \defeq \hspace{-1.em} \underbrace{\B \B \cdots \B}_\text{$\tau$-times composition} \hspace{-1.2em} \bar V, \label{eq:periodic_compensation:operators:optimal} \\
		\begin{split}
		\bar \B_{\bar \eta, \tau}\bar V(z, r) & \defeq \B_{\bar \eta_1}\B_{\bar \eta_2}\cdots\B_{\bar \eta_\tau}\bar V(z, r), \\
		\bar \eta(z, r) & = (\bar \eta_1, \bar \eta_2, \ldots, \bar \eta_\tau).
		\end{split}\label{eq:periodic_compensation:operators:eta}
	\end{align}
\end{subequations}

\begin{theorem}\label{thm:periodic}
	For Problem~\ref{prob:known_dynamics}, let $V_0$ be as in \eqref{eq:known_dynamics:value_iteration},  and let $\B^\tau$ and $\bar \B_{\bar \eta, \tau}$ be as in \eqref{eq:periodic_compensation:operators}.
	If there exists a function $\bar V : \R^{n_z} \times \R^\modelset \to \R$ and positive integer $\tau$ such that
	\begin{equation*}
		\bar V \geq V_0, \quad
		\B^\tau \bar V \leq \bar V,
	\end{equation*}
	then the value iteration $V_0, V_1, \ldots$ converges to a limit $V_\star \leq \bar V$.

	If there exists a $\tau$-periodic control law $\bar \eta$ as in \eqref{eq:periodic_compensation} and \eqref{eq:periodic_compensation:control} such that
	\begin{equation*}
		\bar V \geq V_0, \quad
		\bar \B_{\bar \eta, \tau}\bar V \leq \bar V,
	\end{equation*}
	then $\B^\tau \bar V \leq \bar V$ and the policy $\bar \eta$ achieves a cost, \eqref{eq:known_dynamics:cost}, no greater than
	\begin{equation}\label{eq:periodic:cost}
		\max\{V(z_0, 0), \B_{\bar \eta_\tau}\bar V(z_0, 0), \ldots, \B_{\bar \eta_2}\cdots\B_{\bar \eta_\tau}\bar V(z_0, 0)\},
	\end{equation}
	where $(\bar \eta_1, \bar \eta_2, \ldots, \bar \eta_\tau) = \bar \eta(z_0, 0)$.
\end{theorem}

\begin{corollary}
	If $\bar V$ is the smallest fixed point of $\B^\tau$ greater than $V_0$, then $\bar V = V_\star$.
\end{corollary}

\begin{corollary}
	If one has found a $\bar V \geq V_0$ that satisfies the periodic Bellman inequality, $\B^\tau \bar V \leq \bar V$, then the control law 
	\[
		\bar \eta_{\tau k + s} = \arg\min_u \B_u \B^{\tau - 1 - s} \bar V(z_{\tau k + s}, r_{\tau k + s})
	\]
	achieves a cost, \eqref{eq:known_dynamics:cost}, no greater than
	\[
		\max_{s = 0, 1, \ldots, \tau - 1} \{\B^s \bar V(z_{0}, 0)\}.
	\]
\end{corollary}

\begin{corollary}[1-Periodic]\label{col:periodic_1}
	If there exists a function $\bar V : \R^{n_z} \times \R^\modelset \to \R$ and control law $\bar \eta : \R^{n_z} \times \R^\modelset \to \R^{n_u}$ such that
	$
		\B_{\bar \eta}\bar V \leq \bar V,
	$
	then the value iteration $V_0, V_1, \ldots$ is bounded, and the control policy
	$
		u_t = \bar \eta(z_t, r_t)
	$
	achieves a cost no greater than
	$
		\bar V(z_0, 0)
	$
for Problem~\ref{prob:known_dynamics}.
\end{corollary}

\begin{proof}
	Let $\bar V \geq V_0$ be satisfy $\B^\tau \bar V \leq \bar V$ for some $\tau$.
	As $\B$ is monotone, so is $\B^\tau$.
	Then, for any $k = 0, 1, \ldots$, by monotonicity, $\bar V \geq \B^{k\tau} \bar V \geq \B^{ k\tau} V_0 = V_{k\tau}$.
	By Theorem~\ref{thm:known_dynamics}, the value iteration is monotone, so for any $l = 0, 1, \ldots, \tau - 1$, $\bar V \geq V_{k(\tau + 1)} \geq V_{k\tau + l}$.
	Thus, the value iteration is bounded by $\bar V$.

	For the second part, assume that there exists a $\bar V$ and $\bar \eta$ such that $\B_{\bar \eta, \tau}\bar V \leq \bar V$.
	As $\B_{\bar \eta, \tau}\bar V \geq \B^\tau \bar V$, we have that $\B^\tau \bar V \leq \bar V$.
	It remains to show that the controller $\bar \eta$ achieves a cost no greater than \eqref{eq:periodic:cost}.

	Consider the value iteration starting with $\bar V_0 = \bar V$ and
	\begin{align*}
		\bar V_{(k + 1)\tau}  & = \B_{\bar \eta, \tau}\bar V_k, \\
		\bar V_{k\tau + s + 1} & = \B_{\bar \eta_s}\bar V_{k\tau + s}, \quad s = 0, 1, \ldots, \tau - 2.
	\end{align*}
Fix some $N = 0, 1, \ldots$ and consider the factorization $N = k\tau + s$ where $s = 0, 1, \ldots, \tau - 1$.
We have that
	\begin{align*}
		\sup_{H \in \modelset, d} J^N_{\bar \eta}(z_0, H, d)  & =
	\bar \B_{\eta_1} \bar \B_{\eta_2} \cdots \bar \B_{\eta_s} \bar \B_{\bar \eta, \tau}^k V_0(z_0, 0) \\
								      & \leq
	\bar \B_{\eta_1} \bar \B_{\eta_2} \cdots \bar \B_{\eta_s} \bar \B_{\bar \eta, \tau}^k \bar V(z_0, 0) \\
								      & \leq 
	\bar \B_{\eta_1} \bar \B_{\eta_2} \cdots \bar \B_{\eta_s} \bar V(z_0, 0).
\end{align*}
The bound~\eqref{eq:periodic:cost} follows from taking the supremum over $N = 0, 1, \ldots$ on both sides.
\end{proof}

\subsection{Solution to the Bellman Inequality}\label{sec:explicit:solution}
This section is devoted to an explicit solution to the periodic Bellman inequality, $\B_{\bar \eta, \tau} \bar V \leq \bar V$ in Theorem~\ref{thm:periodic}, in the case of a finite model set
$
	\modelset = \{H_1, \ldots, H_N\}.
$
We parameterize an upper bound of the value function in a set of positive definite matrices $P_{ij} \in \S^{n_z}$,
\begin{subequations}\label{eq:explicit:V_bar}
	\begin{align}
		\bar V(z, r) & = \max_{i,j } \bar V_{ij}(z, r), \\
		\bar V_{ij}(z, r) & = |z|^2_{P_{ij}} + \left(r(H_i) + r(H_j)\right)/2,
	\end{align}
\end{subequations}
where $i, j= 1, \ldots, N$.
We restrict out attention to $\tau$-periodic certainty-equivalence controllers of the form\footnote{If the max is achieved on a set, any selection mechanism will work.}
\begin{equation}\label{eq:explicit:CE_control_law}
	\begin{aligned}
		k(n\tau) & = \argmax_{i \in \{1, \ldots, N\}} r_{n\tau}(H_i), \\
		u_{n\tau + s} & = -K_{k(n\tau)}z_{n\tau + s}, \quad s = 0, \ldots, \tau - 1, \\
	\end{aligned}
\end{equation}
for some matrices $K_1, \ldots, K_N \in \R^{n_u \times n_z}$. 
In the language of Section~\ref{sec:inequalities}, the supervisor, $\bar \eta$, is executed at each $\tau$-th time step and generates the feedback control law to be used over the next $\tau$ time steps:
\[
	\bar \eta(z_{n\tau}, r_{n\tau}) = (\bar \eta_{k(n\tau)}, \bar \eta_{k(n\tau)}, \ldots, \bar \eta_{k(n\tau)}),
\]
where $\bar \eta_k (z, r) = - K_k z$ is the component control law.
\begin{remark}
	The theoretical development in this section does not rely on the gain matrices $K_k$ being constant over each period.
	One could let the supervisor, $\bar \eta$, predetermine a sequence of gain matrices for the next period.
\end{remark}

The Bellman operator acting on a function $\bar V_{ij}$ is the supremum of the quadratic form of the operator $\mathcal G$ acting on the state and the disturbance, $d$:
\begin{multline}
	\B_{-K_kz} \bar V_{ij}(z, r) = \sup_d\Big\{ \quadform_{\mathcal G(P_{ij}, K_k, (H_i + H_j) / 2)}(z, d)\\
	+ (r(H_i) + r(H_j))/2 \Big\},
\end{multline}
where
\begin{multline}\label{eq:G_operator}
\mathcal G(P, K, H) \defeq \bmat{A - BK & G}^\tran P \bmat{A - BK & G} \\
	+ \bmat{I & 0 \\-K & 0 \\ 0 & I}^\tran H \bmat{I & 0 \\-K & 0 \\ 0 & I}.
\end{multline}
We parameterize a bound of the temporal evolution over a period in a sequence of matrices $P_{ij, k}^1, \ldots, P_{ij, k}^\tau \in \S^{n_z}$, so that $\B^s_{\bar \eta_k} V_{ij}(z, r) \leq |z|^2_{P_{ij, k}^s} + (r(H_i) + r(H_j))/2$ for $s = 1, \ldots, \tau$.
This requirement is equivalent to the set of matrix inequalities
\begin{equation}\label{eq:explicit:lmi}
\begin{aligned}
			\bmat{P^1_{ij, k} & 0 \\ 0 & 0} & \succeq \mathcal G(P_{ij}, K_k, (H_i + H_j) / 2), \\
			\bmat{P^{s + 1}_{ij, k} & 0 \\ 0 & 0} & \succeq \mathcal G(P^s_{ij, k}, K_k, (H_i + H_j) / 2).
\end{aligned}
\end{equation}
By the choice of $k$ in~\eqref{eq:explicit:CE_control_law}, we have that $r(H_i) \leq r(H_k)$ for all $i$.
The following theorem formalizes the sufficient condition that  if $P_{ij, k}^\tau \preceq P_{jk}$, then $\bar V$ and $\bar \eta$ fulfills the $\tau$-periodic Bellman inequality.

\begin{theorem}[Explicit solution]\label{thm:explicit}
	For Problem~\ref{prob:known_dynamics} where the model set is finite, $\modelset = \{H_1, \ldots, H_N\}$.
	Assume there exist
	\begin{itemize}
		\item a positive integer $\tau$,
		\item matrices $P_{ij, k}^s = P_{ji, k}^s \in \S^{n_z}$  for $i, j, k = 1, \ldots, N$ and $s = 1, \ldots, \tau$,
		\item positive semidefinite matrices $P_{ij} = P_{ji} \in \S^{n_z}$ for $i, j = 1, \ldots, N$,
		\item gain matrices $K_1, \ldots, K_N \in \R^{n_u \times n_z}$.
	\end{itemize}
	If $P_{ij, k}^\tau \preceq P_{jk}$ and \eqref{eq:explicit:lmi} are fulfilled for all $s = 1, \ldots, \tau - 1$ and $i, j, k = 1, \ldots, N$ except for $i \neq j = k$, then the value approximation $\bar V$ in \eqref{eq:explicit:V_bar} with the certainty-equivalence control law $\bar \eta$ in \eqref{eq:explicit:CE_control_law} fulfills the periodic Bellman inequality, $\B_{\bar \eta, \tau} \bar V \leq \bar V$.
Furthermore, the $\bar \eta$ achieves an objective value no greater than
\begin{multline*}
	\max\Big\{z_0^\tran Z_{ij} z_0 + (r(H_i) + r(H_j))/2, 
	Z_{ij} \in \{P_{ij}\} \cup \{P_{ij, k}^s\},\\
	i, j, k \in \{1, \ldots, N\}, s \in \{1, \ldots, \tau-1\Big\}
\end{multline*}
\end{theorem}
\begin{proof}
	As $P_{ij} \succeq 0$, we have $\bar V \geq V_0$.
	The inequalities \eqref{eq:explicit:lmi} imply that $\B_{\bar \eta_k}^s V_{ij} \leq |z|^2_{P_{ij, k}^s} + (r(H_i) + r(H_j))/2$ for all $s = 1, \ldots, \tau$.
	By the choice of $k$ in~\eqref{eq:explicit:CE_control_law}, we have that $r(H_i) \leq r(H_k)$ for all $i$ and
	$\B_{\bar \eta, \tau} \bar V_{ij} \leq \max_{k} (\bar V_{jk} + (r(H_i) + r(H_j))/2) \leq \bar V$.
	Finally as $\max_{ij}\B_{\bar \eta, \tau} \bar V_{ij} = \B_{\bar \eta, \tau} \max_{ij}\bar V_{ij}$, the conclusions follow from Theorem~\ref{thm:periodic}.
\end{proof}

\begin{remark}
	Theorem 3 in~\cite{Rantzer2021L4DC} is obtained as a corollary of Theorem~\ref{thm:explicit} by substituting $\hat A = 0$, $\hat B = 0$ and $\hat G = I$ and $H_M$ from \eqref{eq:state_feedback:H} into $\mathcal G$, taking $\tau = 1$ and replacing \eqref{eq:explicit:lmi} with their lower Schur complements.
\end{remark}

\begin{remark}\label{rem:explicit}
	The inequalities \eqref{eq:explicit:lmi} are affine in $P_{ij}$ and $P_{ij, k}^s$, but are not convex in $K_i$.
	One heuristic approach to solve the inequalities is to first solve the linear-quadratic problem associated with each model $i$ to obtain $K_i$,
\[
	\inf_\mu\sup_{d, N}\left\{\sum_{t = 0}^{N-1}\quadform_{H_i}(z_t, u_t, d_t) \right\}.
\]
	Then use standard optimization software for semidefinite programming to search for $P_{ij}$ and $P_{ij, k}^s$, holding $K_i$ fixed.
	This approach was suggested in~\cite{Rantzer2021L4DC} and is also used in the examples in Section~\ref{sec:examples}.
\end{remark}

\section{Examples}\label{sec:examples}
\subsection{State-Feedback, Delays and Periodic Compensation}\label{sec:examples:statefeedback}
This section studies state-feedback control of the delayed discrete-time integrator where the sign of the gain is unknown.
The dynamics can be modeled in two ways, either the sign uncertainty is incorporated into the state matrix or the input:
\begin{align}
			x_{t + 1} & = \underbrace{\bmat{1 & i \\ 0 & 0}}_{A_i}x_t + \underbrace{\bmat{0 \\ 1}}_{B}u_t + w_t,\quad i = \pm 1\label{eq:sf:linear:delay:state} \\
			x_{t + 1} & = \underbrace{\bmat{1 & 1 \\ 0 & 0}}_{A}x_t + \underbrace{\bmat{0 \\ i}}_{B_i}u_t + w_t, \quad i = \pm 1, \label{eq:sf:linear:delay:input}
\end{align}
with $Q = I$ and $R = I$.
Although the input to output ($x_1$) behavior of the systems~\eqref{eq:sf:linear:delay:state} and~\eqref{eq:sf:linear:delay:input} are identical, from a control perspective, they are significantly different.
The difference lies in that an impulse in the controlled input at time $t$ will reveal information about the sign of $B_i$ in~\eqref{eq:sf:linear:delay:input} at time $t + 1$ but information about $A_i$ in~\eqref{eq:sf:linear:delay:state} not until $t + 2$.
This difference is reflected in the smallest period, $\tau$, for which the periodic Bellman inequality can be satisfied.

We computed $K_1$ and $K_2$ according to Remark~\ref{rem:explicit} and solved the conditions in Theorem~\ref{thm:explicit} using \texttt{MOSEK}.
We find that for \eqref{eq:sf:linear:delay:input}, the conditions are satisfied for $\tau = 1$ and $\gamma = 6$.
Our software implementation cannot find a solution for $\tau = 1$ for the system in~\eqref{eq:sf:linear:delay:state} and it is not until $\tau = 2$ that a solution, with $\gamma = 11.2$, is found.
\begin{prop}
	Given $A_i$ and $B$ in \eqref{eq:sf:linear:delay:state}, $\gamma > 0$, $Q \in \R^{2 \times 2}$ be positive definite and $R > 0$.
	Then there does not exist matrices $K_i \in \R^{1 \times 2}$ and positive definite $P \in \R^{2 \times 2}$ such that $A_i - BK_i$ are Schur stable and,
	\begin{equation}\label{eq:example:propineq}
		\bmat{P & 0\\ 0 & 0} \geq \mathcal G(P, K_k, (H_1 + H_2) / 2)
	\end{equation}
	for both $(i, j, k) = (1, -1, 1)$, and $(-1, 1, -1)$.
\end{prop}
\begin{proof}	
	Taking the Schur complement of \eqref{eq:example:propineq} we get the equivalent condition that for all $x$
	\begin{multline*}
		|x|_P^2 \geq |x|_Q^2 + |K_kx|^2_R - \gamma^2 |(A_i - B_iK_k - A_j + B_jK_k)x / 2|^2 \\
		+ |(A_i - B_iK_k + A_j - B_jK_k)x / 2|^2_{(P^{-1} - \gamma^{-2}I)^{-1}}.
	\end{multline*}
	Note that $(P^{-1} - \gamma^{-2}I)^{-1} \succ P$.
	Let $K_l = \bmat{k_l^1 & k_l^2}$.
	By Lemma~\ref{lemma:stability_conditions} in the Appendix, that $lk^1_l > 0$.
	Furthermore, we have for $i \neq j$
	\begin{equation*}
		A_i - B_i K_l + A_j - B_j K_l = \bmat{
			1 & 0 \\
			-k_l^1 & - k_l^2
		}.
	\end{equation*}
	Thus,
	\begin{multline*}
		(A_i - B_i K_l + A_j - B_j K_l)^\tran P (A_i - B_i K_l + A_j - B_j K_l) \\
		 = \bmat{
			 p_{11} - 2k_l^1 p_{21} + (k_l^1)^2p_{22} & * \\ * & *
		},
	\end{multline*}
	where $p_{12} = p_{21}$ and
	\[
		P = \bmat{p_{11} & p_{12} \\ p_{21} & p_{22}}.
	\]
	We also note that
	\[
		- \gamma^2 |(A_i - B_i K_l - A_j + B_j K_l)x / 2|^2 = - \gamma^2 x_2^2.
	\]
	For $x = \bmat{x_1 & 0}^\tran$, we get
	\begin{align*}
		& |x|_P^2 - |x|_Q^2 + |K_kx|^2_R 
		- \gamma^2 |(A_i - B_i K_l - A_j + B_j K_l)x / 2|^2 \\
		& \quad +|(A_i - B_i K_l + A_j - B_j K_l) x / 2|^2_{(P^{-1} - \gamma^{-2}I)^{-1}} \\
		&\leq |x|_P^2 - |x|_Q^2 + |K_kx|^2_R 
		- \gamma^2 |(A_i - B_i K_l - A_j + B_j K_l)x / 2|^2 \\
		& \quad +|(A_i - B_i K_l + A_j - B_j K_l) x / 2|^2_P \\
		& = 2k_l^1 p_{21} - (k_l^1)^2 p_{22} - q_{11} - (k_l^1)^2 R
	\end{align*}
	As $(k_l^1)^2 p_{22} + q_{11} + (k_l^1)^2R > 0$ and as $k_1^1$ and $k_{-1}^1$ have opposite signs, we have that the last line is smaller than zero for $l = 1$, $l = -1$ or both.
\end{proof}

	\begin{figure*}
		\centering
		\forJournal{
			\includegraphics{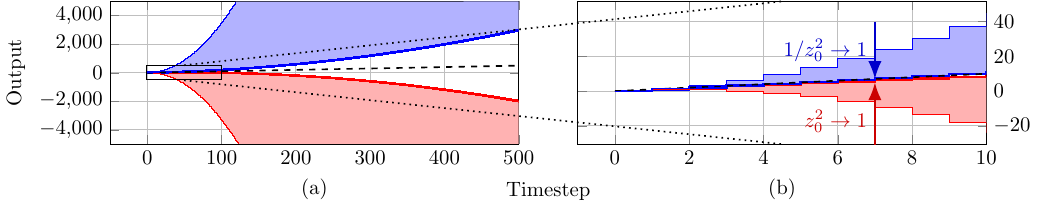}
		}
		\forArxiv{
			\includegraphics{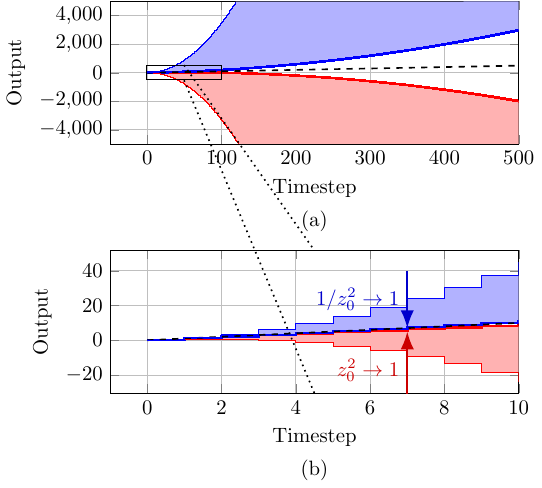}
		}
		\caption{Step responses for the systems $G_\text{mp}$ (blue) and $G_\text{nmp}$ (red) as $z_0$ varies continuously from $1.4$ to $1.01$ (thick line). The dashed line $y = t$ corresponds to an integrator. As $z_0 \to 1$, note that the short-term behavior of both realizations resemble that of the integrator. The asymptotic behavior, however, do not.
		\label{fig:step_responses}
		}
	\end{figure*}
	\begin{figure}
		\centering
		\includegraphics{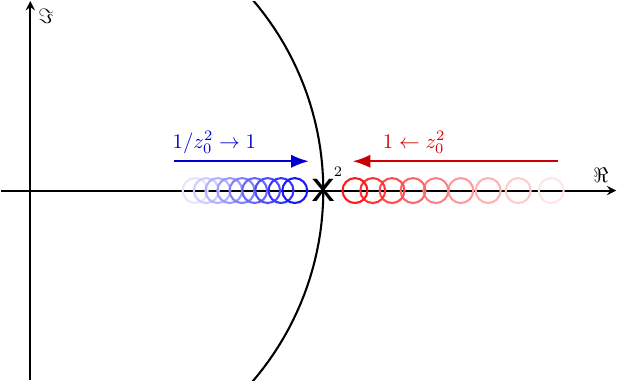}
		\caption{Pole-zero map of the systems $G_\text{mp}$ (blue) and $G_\text{nmp}$ (red) as $z_0$ approaches $1.0$. The double pole at $z = 1$ is unstable, and the zero at $z = 1/z_0^2$ (blue) is minimum phase. The zero at $z = z_0^2$ (red) is nonminimum phase. The zeros are reflections in the unit circle.}
		\label{fig:pzmap}
	\end{figure}
	\begin{figure}
		\centering
		\includegraphics{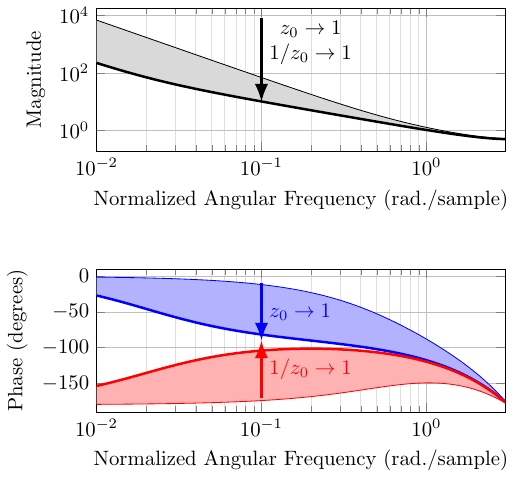}
		\caption{Bode plot of the systems $G_\text{mp}$ (blue) and $G_\text{nmp}$ (red) as $z_0$ varies continuously from $1.4$ to $1.01$ (thick).
		The systems' magnitude responses are equal for a fixed $z_0$, but they differ in phase. This difference is negligable for high frequncies.}
		\label{fig:bode_plot}
	\end{figure}
	\begin{figure*}
		\centering
		\forJournal{
		\includegraphics{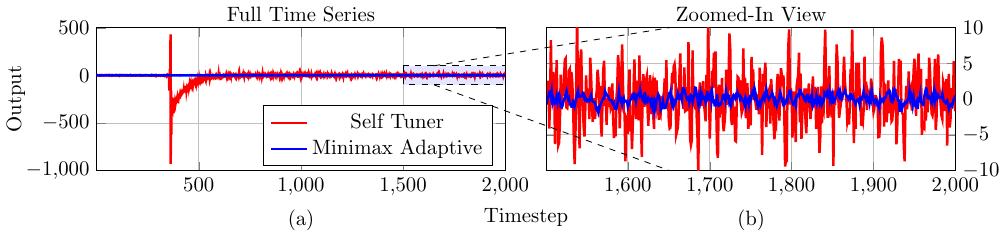}
		}
		\forArxiv{
			\includegraphics{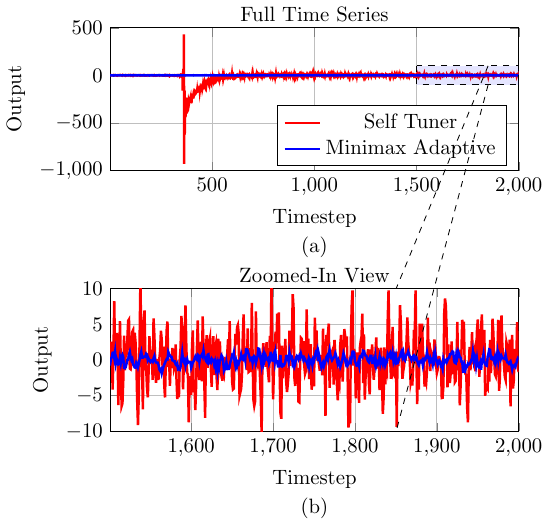}
		}
		\caption{Simulations of the double integrator with uncertain approximate pole cancellation. The system is in feedback with the periodic certainty-equivalence controller of Section~\ref{sec:explicit} (Minimax Adaptive, blue) and the self-tuning LQG controller described in~\cite[Chapter 4]{Astrom1995Adaptive} (Self tuner, red). The output is shown the entire duration in (a), and a shorter snapshot in (b). Note the spike for the self-tuning regulator at time-step 355. These spikes do not occur with the periodic controller.}
		\label{fig:mac_str_signals}
	\end{figure*}
	\begin{figure}
		\centering
		\includegraphics{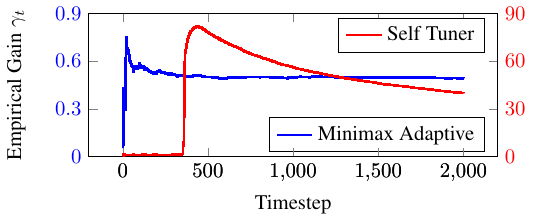}
		\caption{Evolution of the empirical gain $\gamma_t = \sum_{\tau = 0}^t\left(|x_\tau|_Q^2 + |u_\tau|_R^2 \right) / \left(|x_0 - \hat x_0|^2_S +\gamma^2 \sum_{\tau=0}^t |(w_\tau, v_\tau)|^2\right) $ for the double integrator with uncertain approximate pole cancellation. The Minimax Adaptive controller (blue, left scale) initially has a higher $\gamma_t$ than the self-tuning LQG controller (red, right scale), but due to the build up in the second-order mode, the self tuner spikes at $t = 355$.}
		\label{fig:mac_str_gammas}
	\end{figure}
	\begin{figure}
		\centering
		\includegraphics{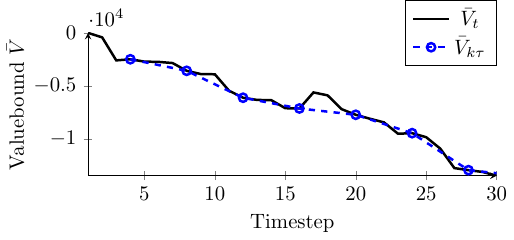}
		\caption{Evolution of the value function approximation for the double integrator with uncertain approximate pole cancellation. The value function is shown for the Minimax Adaptive controller and the periodic decrease, with period $\tau = 4$, is highlighted with circles.}
		\label{fig:evilV}
	\end{figure}

\subsection{Approximate Unstable Pole Cancellation}\label{sec:examples:output}
	We conclude the examples by synthesizing a controller for the double integrator with uncertain approximate pole cancellation.
	See the pole-zero map in Fig.~\ref{fig:pzmap}.
	This corresponds to an approximate cancellation of the unstable pole.
	The step responses of the system in Fig.~\ref{fig:step_responses}, and Bode plots in Fig.~\ref{fig:bode_plot} indicate that the high-frequency behavior of $G_\text{mp}$ and $G_\text{nmp}$ are similar to an integrator, but that the low-frequency asymptotes are different.

	The minimum phase system, $G_\text{mp}$, and the nonminimum-phase system, $G_\text{nmp}$, have state-space realizations $(A, B, C_\text{mp}, D, G)$ and $(A, B, C_\text{nmp}, D, G)$ respectively, where
	\begin{align*}
		A & = \bmat{1 & 1\\ 0 & 1}, & B & = \bmat{0 \\ 1}, \\
		C_\text{mp} & = \bmat{-1/z_0 + z_0 \\ z_0}^\tran, & C_\text{nmp} & = \bmat{-z_0 + 1/z_0 \\ 1/z_0}^\tran.
	\end{align*}
	Here $z_0 = 1.01$.
	The LMIs \eqref{eq:explicit:lmi} have solutions for $\tau = 4$ and $\gamma = 20$ with
	$G = I / 100$, $Q = I / 100$, $R = I / 100$ and $D = 1 / 10$.
	The matrices were scaled down for numerical stability in the optimization.
	We evaluate the performance of the periodic certainty-equivalence controller and compare to self-tuning LQG controller described in~\cite[Chapter 4]{Astrom1995Adaptive} by simulating the nonminimum-phase system with $w_t$ and $v_t$ normally distributed with zero mean and unit variance.

	Time series of the output signal are shown in Fig.~\ref{fig:mac_str_signals}. 
	The self-tuning controller has a spike at time-step $355$ which is due to its inability to act against the growth of the mode associated with the second integrator before it starts dominating the output.
	The minimax adaptive controller does not lead to such spikes.
	We also show the evolution of $\gamma_t$ in Fig.~\ref{fig:mac_str_gammas} and the evolution of the value function in Fig.~\ref{fig:evilV}.

\section{Conclusions}
We conclude with a few words about the limitations of this work and promising research directions.
\begin{enumerate}
	\item Including pathological hypotheses in the model set, such as $G_\text{nmp}$ in~\eqref{eq:intro:nmp} severely impacts the performance guarantee, even when the actual realization of the system is well-behaved.
		The actual performance may be much better, as indicated by comparing the empirical gain in Fig.~\ref{fig:mac_str_gammas}, $\gamma_t < 0.9$, with the guarantee $\gamma < 20$, and we do not yet fully understand the effects of including pathological hypotheses.
		Integrating the framework of Goel \textit{et al.}~\cite{Goel2023Competitive} with the results of this article seems promising to address the pathological hypotheses in the model set. 
		As the authors reformulate regret optimization and competitive ratio optimization as $\Hinf$ synthesis problems, our results could be integrated to compute suboptimal control policies in the case of parametric uncertainty.
	\item We assumed that the model set was finite.
		Even though finite sets can approximate compact sets of models, our results do not inform how to choose the approximate models and how to quantify the approximation error.
	Theorem~\ref{thm:output_feedback:reduction} shows that by constructing one $\Hinf$-observer for each model in the model set, the output feedback problem can be exactly reduced to an instance of Problem~\ref{prob:known_dynamics}.
		If the model set is infinite but compact, one could instead construct a robust $\Hinf$-observer for each element in a finite cover and approximate the output feedback problem with an instance of Problem~\ref{prob:known_dynamics}.
	\item Section~\ref{sec:examples:statefeedback} demonstrates that when delays are present, the value function approximation \eqref{eq:explicit:V_bar} does not capture the probing effect of the control policy. 
		The introduction of periodicity in the control policy, as in Section~\ref{sec:explicit}, mitigates this problem as it allows for information gathering over a longer time horizon.
		Capturing the probing effect of the control policy in the value approximation is crucial for obtaining tighter performance bounds. 
		Numerical studies of the value iteration could provide insight into this.
		One could also investigate using reinforcement learning to approximate the value function.
\end{enumerate}

\appendix

\begin{lemma}\label{lemma:stability_conditions}
        The system
	\begin{equation}
		x_{t+1} = \bmat{1 & i \\ 0 & 1} x_t + \bmat{0 \\ 1} u_t,
	\end{equation}
	where $i \in \pm 1$ and 
        \[  
                u[t] = K_i x[t]
        \]
        is asymptotically stable if, and only if
        \begin{subequations} \label{eq:stability_conditions}
                \begin{equation}\label{eq:stability_conditions_k1} 
                        ik^i_1 > 0,
                \end{equation}
                \text{and}
                \begin{equation}\label{eq:stability_conditions_k2} 
                        ik^i_1 - 1 < k^i_2 < 1 + \frac{ik^i_1}{2}.
                \end{equation}
        \end{subequations}
\end{lemma}
\begin{proof}
        The characteristic polynomial of the closed-loop system is
        \[  
	(z - 1)(z + k^i_2) + ik^i_1 = z^2 + (k^i_2 - 1)z + ik^i_1 - k^i_2.
        \]
        The inequalities \eqref{eq:stability_conditions} follow by the Jury stability criterion.
\end{proof}

\section*{Acknowledgment}
The authors would like to thank their colleague Venkatraman Renganathan for valuable feedback on the manuscript.
ChatGPT~\cite{Openai2024Chatgpt} was used for light editing of this text.
\forJournal{
	\bibliographystyle{IEEEtranDOI}
}
\forArxiv{
	\bibliographystyle{IEEEtranDOI}
}
\bibliography{IEEEabrv, references}             
\forJournal{
	\begin{IEEEbiography}[{\includegraphics[width=1in,height=1.25in,clip,keepaspectratio]{pics/olle_small.jpg}}]{Olle Kjellqvist} (Student Member, IEEE) %
was born in Halmstad, Sweden, in 1990. 
He received the M.Sc degree in engineering physics from Lund University, Lund, Sweden, in 2018, and the Tekn. Lic. degree from Lund University in 2022.

From 2018 to 2019 he was with Accenture, Copenhagen, Denmark, as a software consultant.
He is currently pursuing the Ph.D. degree in automatic control at Lund University.
His research interests include adaptive control, robust control, learning-based control, and control of large-scale systems.

Lic. Kjellqvist is a member of Engineers of Sweden and the ELLIIT Strategic Research Area at Lund University.
\end{IEEEbiography}%
\begin{IEEEbiography}[{\includegraphics[width=1in,height=1.25in,clip,keepaspectratio]{pics/rantzer_2_big.png}}]{Anders Rantzer} (Fellow, IEEE) %
	Anders Rantzer was appointed professor of Automatic Control at Lund
University, Sweden, after a PhD at KTH Stockholm in 1991 and a postdoc
1992/93 at IMA, University of Minnesota. The academic year of 2004/05 he
was visiting associate faculty member at Caltech and 2015/16 he was
Taylor Family Distinguished Visiting Professor at University of
Minnesota. Rantzer has served as chairman of the Swedish Scientific
Council for Natural and Engineering Sciences as well as the Royal
Physiographic Society of Lund. He is a Fellow of IEEE and member of the
Royal Swedish Academy of Engineering Sciences. His research interests
are in modeling, analysis and synthesis of control systems, with
particular attention to scalability, adaptation and applications in
energy networks.
\end{IEEEbiography}

}
\end{document}